\documentclass[a4paper,12pt]{amsart}
\usepackage{ifthen}
\usepackage{mathrsfs}
\numberwithin{equation}{section}
\nonstopmode
\newtheorem{thm}{Theorem}[section]

\newtheorem{lem}[thm]{Lemma}
\newtheorem{prop}[thm]{Proposition}

\theoremstyle{definition}

\newtheorem{example}[thm]{Example}

\newenvironment{pf}[1][]{%
 \vskip 3mm
 \noindent
 \ifthenelse{\equal{#1}{}}%
  {{\slshape Proof. }}%
  {{\slshape #1.} }%
 }%
{\qed\bigskip}

\newcounter{alphabet}
\newcounter{tmp}

\newcommand{\A}{{\mathcal A}}

\newcommand{\C}{{\mathbb C}}
\newcommand{\D}{{\mathbb D}}

\newcommand{\M}{{\mathcal M}}

\newcommand{\N}{{\mathcal N}}

\newcommand{\R}{{\mathbb R}}

\newcommand{\es}{{\mathcal S}}

\newcommand{\ZF}{{\mathcal{ZF}}}

\newcommand{\sphere}{{\widehat{\mathbb C}}}

\renewcommand{\Im}{\,{\operatorname{Im}\,}}

\renewcommand{\Re}{{\operatorname{Re}\,}}

\newcommand{\Int}{{\operatorname{Int}\,}}

\newcommand{\inv}{^{-1}}

\newcommand{\arctanh}{{\operatorname{arctanh}}}

\renewcommand{\arg}{\,{\operatorname{arg}\,}}

\newcommand{\Log}{{\operatorname{Log}\,}}

\newcounter{minutes}
\divide\time by 60
\newcounter{hours}
\multiply\time by 60
\addtocounter{minutes}{-\time}

\begin{document}
\bibliographystyle{amsplain}
\title{
On univalence of the power deformation $z(f(z)/z)^c$
}


\author[Y.~C.~Kim]{Yong Chan Kim}
\address{Department of Mathematics Education, Yeungnam University, 214-1 Daedong
Gyongsan 712-749, Korea}
\email{kimyc@ynu.ac.kr}
\author[T. Sugawa]{Toshiyuki Sugawa}
\address{Division of Mathematics, Graduate School of Information Sciences,
Tohoku University, Aoba-ku, Sendai 980-8579, Japan}
\email{sugawa@math.is.tohoku.ac.jp}
\begin{abstract}
In this note,
we mainly concern the set $U_f$ of $c\in\C$ such that the power deformation
$z(f(z)/z)^c$ is univalent in the unit disk $|z|<1$ for a given
analytic univalent function $f(z)=z+a_2z^2+\cdots$ in the unit disk.
We will show that $U_f$ is a compact, polynomially convex subset of the
complex plane $\C$ unless $f$ is the identity function.
In particular, the interior of $U_f$ is simply connected.
This fact enables us to apply various versions of the $\lambda$-lemma
for the holomorphic family $z(f(z)/z)^c$ of injections parametrized
over the interior of $U_f.$
We also give necessary or sufficient conditions for $U_f$ 
to contain $0$ or $1$ as an interior point.
\end{abstract}
\thanks{
The first author was supported by Yeungnam University (2011).
The second author was supported in part by JSPS Grant-in-Aid for 
Scientific Research (B) 22340025.
}
\maketitle

\section{Introduction}

Let $\A$ be the class of analytic functions on the unit disk
$\D=\{z\in\C: |z|<1\}.$
We denote by $\A_0$ its subclass consisting of functions $h$
normalized by $h(0)=1.$
The set $\A_0^\times$ of invertible elements in $\A_0$ with respect to
pointwise multiplication is nothing but the set of non-vanishing functions
in $\A_0.$
For $h\in\A_0^\times,$ define $\Log h$ to be the analytic branch
of $\log h$ on $\D$ determined by the condition $\Log h(0)=0.$
The set of functions $f$ in $\A$ with the representation $f(z)=zh(z)$
for some $h\in\A_0$ (resp.~$h\in\A_0^\times$) will be designated by
$\A_1$ (resp.~$\ZF$).
In other words, $f\in\A$ belongs to $\A_1$ if and only if $f(0)=0, f'(0)=1;$
and $f\in\A$ belongs to $\ZF$ if and only if $f\in\A_1$ and
$f(z)\ne0$ for $0<|z|<1.$
We denote by $\es$ the set of univalent functions in $\A_1.$
Note that $\es$ is contained in $\ZF.$

In \cite{KS11PT}, the authors investigated the power deformation
$$
K_c[f](z)=z\left(\frac{f(z)}{z}\right)^c
$$
for $f\in\ZF$ and $c\in\C.$
Here and in what follows, the power $h^c$ will be defined as
$\exp(c\, \Log h)$ for $h\in\A_0^\times$ and $c\in\C.$
We determined the sets $[\M,\N]_K=\{c\in\C: K_c[\N]\subset\M\}$
for various subclasses $\M,\N$ of $\es$ in \cite{KS11PT}.
In the present note, we focus our attention on the set
$$
U_f=[\{f\},\es]_K=\{c\in\C: K_c[f]~\text{is univalent on }\D\}
$$
for $f\in\ZF.$
For instance, for the Koebe function $\kappa(z)=z/(1-z)^2,$ we have
$U_\kappa=\{c: |c-1/2|\le 1/2\}$ (see the remark right after the proof
of Theorem 1.1 in \cite{KS11PT}).
By the property $K_c\circ K_{c'}=K_{cc'},$ we have the relation
$U_{K_c[f]}=c\inv\cdot U_f$ (see \cite[Lemma 2.1]{KS11PT}).

Note that $U_f$ has an interior point when $f$ is a starlike univalent
function or, more generally, a spirallike function (see \cite{KS11PT}).
We are motivated, in part, by the fact that the interior $\Int U_f$ serves as
a parameter region of the holomorphic family $K_c[f]$ of injections on $\D.$
Therefore, we could relate the present study to the theory of quasiconformal
mappings and Teichm\"uller spaces.

We recall here the notion of holomorphic motions.
A holomorphic motion of a subset $E$ of the Riemann sphere 
$\sphere=\C\cup\{\infty\}$ over a domain $D$ with base point $c_0$
is a map $F:D\times E\to\sphere$ with the following three properties:
\begin{enumerate}
\item
$F(c,\cdot): E\to \sphere$ is injective for each $c\in D,$
\item
$F(\cdot, z): D\to\sphere$ is holomorphic for each $z\in E,$ and
\item
$F(c_0,z)=z$ for $z\in E.$
\end{enumerate}
This simple notion appeared only recently
in a paper \cite{MSS83} by Ma\~n\'e, Sad and Sullivan
to study complex dynamics, and afterwards, it found
many applications in various branches of complex analysis.
We summarize necessary results concerning holomorphic motions
in Section \ref{sec:qc}.

In the present note, we will show the following.

\begin{thm}\label{thm:Uf}
Suppose that $f\in\ZF$ is not the identity function.
Then $U_f$ is a compact, polynomially convex set in $\C$ with $0\in U_f.$
\end{thm}

Note that $U_f=\C$ when $f$ is the identity function.
We recall here that a compact set $E$ in $\C$ is polynomially convex
if and only if $\C\setminus E$ is connected 
(see \cite[Chap.~VII, Prop.~5.3]{Conway:FA} for instance).
The latter condition is also known as a characterization of the Runge property
in dimension one.
In particular, we see that each connected component of the interior
$\Int U_f$ is simply connected.

\begin{thm}\label{thm:qc}
Let $D$ be a connected component of $\Int U_f$ for a non-identity
function $f\in\ZF.$
Then the family of univalent functions $f_c=K_c[f]$
over $c\in D$ is quasiconformally homogeneous in $\C:$ more precisely,
for each pair of points $c_0$ and $c_1$ in $D$ there exists a 
$\tanh d_D(c_0,c_1)$-quasiconformal
automorphism $g$ of $\C$ such that $f_{c_1}=g\circ f_{c_0}$ on $\D.$
\end{thm}

Here, $d_D$ denotes the hyperbolic distance in $D$ induced by
the hyperbolic metric of constant curvature $-4.$
For instance, $d_\D(0,z)=\frac12\log\frac{1+|z|}{1-|z|}=\arctanh|z|.$
A mapping $g:D_1\to D_2$ between domains $D_1$ and $D_2$ in $\sphere$
is called {\it $k$-quasiconformal} if $g$ is a homeomorphism
with locally integrable partial derivatives on 
$D_1\setminus\{\infty,g\inv(\infty)\}$ such that 
$|\partial_{\bar z}g|\le k|\partial_zg|$ a.e.~in $D_1$
for a constant $0\le k<1.$

The quasiconformal homogeneity in $\C$ implies, for instance, that
$f_{c_1}$ is bounded on $\D$ precisely when so is $f_{c_0}$ for
$c_0,c_1\in D.$

A key step to prove the last theorem is the fact that
$F(c,z)=f_c(z)$ satisfies conditions (1) and (2) in the definition
of holomorphic motions.
We note that condition (3) is also satisfied when $0\in D\subset\Int U_f.$
However, there is no guarantee that the holomorphic family $f_c$ of injections
over $D$ contains $f$ itself.
This happens when $1\in D\subset\Int U_f.$
If $\{0,1\}\subset D,$ then $f$ is a quasiconformal deformation of
the identity mapping, and therefore $f$ extends to a quasiconformal
automorphism of $\C.$
Here we have conditions for these situations.

\begin{thm}\label{thm:char}
Let $f$ be a non-identity function in $\ZF.$
\begin{enumerate}
\item
$0\in\Int U_f$ if and only if the function $zf'(z)/f(z)$ is bounded on $\D.$
\item
Suppose that $1\in\Int U_f.$
Then $f\in\es$ and the function $zf'(z)/f(z)$ is bounded away from $0$ on $\D.$
\end{enumerate}
\end{thm}

The converse is not true in general in the second
assertion of the last theorem (see Lemma \ref{lem:LU} and
Example \ref{ex:1} below).
Though we do not find so far a sufficient condition general enough,
we have several geometric conditions for $f$ to have the property
$1\in\Int U_f.$
For instance, it suffices to assume that $f$ is starlike of order $\alpha$
for some $\alpha>0.$
See \cite{KS11PT} for details.

We note here that $\Int U_f$ might be empty.
On the other hand, $\Int U_f$ may have many components.
We will show the following result.

\begin{thm}\label{thm:two}
There does exist a function $f\in\es$ such that $\Int U_f$ consists
of at least two connected components.
\end{thm}

We briefly describe the organization of the present note.
In Section 2, we prove Theorems \ref{thm:Uf} and \ref{thm:char}.
There, key ingredients are an idea of \v Zuravlev \cite{Zur80}
and a fundamental relation in \eqref{eq:fund} (see also \cite{KS11PT})
between a set $LU_f$ containing $U_f$ and the variability region $V(f)$
of $zf'(z)/f(z).$
Section 3 is a short section giving a version of the $\lambda$-lemma
and a proof of Theorem \ref{thm:qc}.
In Section 4, we prove Theorem \ref{thm:two} and give a couple of related
results.
To prove the theorem, we prepare a univalence criterion (Lemma \ref{lem:univ}),
which may be of independent interest.

\section{Proof of Theorems \ref{thm:Uf} and \ref{thm:char}}

Univalence is a global property of a function so that
it is not easy to check.
Therefore, it is helpful to consider local univalence instead
as in \cite{KS11PT}.
Recall that an analytic function $f$ is locally univalent at $z_0$
if and only if $f'(z_0)\ne0.$
For a function $f$ in $\ZF,$ we set
$$
LU_f=\{c\in\C: K_c[f]~\text{is locally univalent on }\D\}.
$$
Obviously, $U_f\subset LU_f.$

We now set $f_c=K_c[f]$ for brevity.
A simple computation gives us the relation
\begin{equation}\label{eq:fc}
\frac{zf_c'(z)}{f_c(z)}=1-c+c\frac{zf'(z)}{f(z)}.
\end{equation}
Hence, for a point $z_0\in\D,$ $f_c'(z_0)=0$ if and only if
$z_0f'(z_0)/f(z_0)=(c-1)/c,$ equivalently,
$c=T(z_0f'(z_0)/f(z_0)),$ where
\begin{equation}\label{eq:T}
T(w)=\frac1{1-w}.
\end{equation}
In this way, we have the fundamental relation
\begin{equation}\label{eq:fund}
LU_f=\C\setminus T(V(f)),
\end{equation}
where $V(f)$ is the image of $\D$ under the function $zf'(z)/f(z);$ namely
$$
V(f)=\{zf'(z)/f(z): z\in\D\}.
$$

We need to recall the Grunsky theorem to prove polynomial convexity of $U_f.$
The reader may refer to \cite{Pom:univ} for details.
The Grunsky coefficients $b_{jk}$ of $f\in\A_1$ are defined by expansion
in the form
$$
\log\frac{1/f(z)-1/f(w)}{1/z-1/w}
=-\sum_{j,k=1}^\infty b_{jk}z^jw^k
$$
of double power series convergent in $|z|<\delta, |w|<\delta$
for small enough $\delta>0.$
Indeed, we can take $\rho$ as $\delta$ when $f(z)$ is univalent on the
disk $|z|<\rho.$
The Grunsky theorem says that $f$ is univalent on $\D$ if and only if
$$
\left|\sum_{j,k=1}^N b_{jk}x_jx_k\right|
\le\sum_{j=1}^N\frac{|x_j|^2}{j}
$$
for any positive integer $N$ and any vector $(x_1,\dots,x_N)\in\C^N.$

We are now ready to prove Theorems \ref{thm:Uf} and \ref{thm:char}.

\begin{pf}[Proof of Theorem \ref{thm:Uf}]

Let $f\in\ZF$ be a non-identity function.
Then $zf'(z)/f(z)$ is not constant (and thus an open mapping).
Therefore, $V(f)$ is an open neighbourhood of $1,$ which
implies that $T(V(f))$ is an open neighbourhood of $\infty.$
Now the relation \eqref{eq:fund} yields that $LU_f$
is a compact subset of $\C.$
Since $U_f\subset LU_f,$ we conclude that $U_f$ is bounded.
The Hurwitz theorem implies that $U_f$ is closed.
Hence, $U_f$ is compact.

We next show that $U_f$ is polynomially convex by employing
the idea of \v Zuravlev \cite{Zur80}.
Suppose, to the contrary, that $\C\setminus U_f$ has a bounded
component $\Delta.$
Then we note that $\partial\Delta\subset U_f$ and $\Delta\cap U_f=\emptyset.$
We denote by $b_{jk}(c)$ the Grunsky coefficients of
the function $f_c=K_c[f].$
Then $b_{jk}(c)$ is a holomorphic function in $c$ for each
pair of $j$ and $k.$
(Indeed, it is not difficult to see that $b_{jk}(c)$ is a polynomial in $c.$)
By the Grunsky theorem, for each $(x_1,\dots,x_N)\in\C^N,$ the inequality
\begin{equation}\label{eq:G}
\left|\sum_{j,k=1}^N b_{jk}(c)x_jx_k\right|
\le\sum_{j=1}^N\frac{|x_j|^2}{j}
\end{equation}
holds for $c\in\partial\Delta\subset U_f.$
By the maximum modulus principle for analytic functions,
we see that the inequality \eqref{eq:G} still holds for all
$c\in\Delta.$
Therefore, by the converse part of the Grunsky theorem,
we conclude that $f_c$ is univalent for $c\in\Delta.$
This means that $\Delta\subset U_f,$ which is a contradiction.

The assertion $0\in U_f$ is trivial.
The proof is now complete.
\end{pf}

\begin{pf}[Proof of Theorem \ref{thm:char}]
Let $f\in\ZF$ be a non-identity function.
First we assume that $0\in\Int U_f.$
Then, $0\in\Int LU_f,$ which implies that $0$ is an exterior point
of $T(V(f)).$
Since $T(\infty)=0,$ we have that $V(f)$ is bounded.

Conversely, we assume that $zf'(z)/f(z)$ is bounded.
Then, by \eqref{eq:fc}, the range of $zf_c'(z)/f_c(z)$ shrinks
to the point $1$ when $c$ approaches $0.$
In particular, $\Re[zf_c'(z)/f_c(z)]>0$ for sufficiently small $c.$
In this case, $f_c$ is a starlike univalent function.
Therefore, a neighbourhood of $0$ is contained in $U_f.$
Thus the first part of the theorem is confirmed.

Finally, we assume that $1\in\Int U_f.$
Then $1\in U_f;$ namely, $f$ is univalent.
Since $U_f\subset LU_f$ it is enough to show the following lemma
to finish the proof.
\end{pf}

\begin{lem}\label{lem:LU}
For a function $f\in\ZF,$
$1\in\Int LU_f$ if and only if $f(z)/zf'(z)$ is bounded.
\end{lem}

\begin{pf}
In view of the relation \eqref{eq:fund}, we see that
$1\in\Int LU_f$ precisely when $1$ is an exterior point of $T(V(f)).$
Since $T(0)=1,$ the last condition means that $0$ is an exterior point of $V(f);$
namely, $zf'(z)/f(z)$ is bounded away from $0.$
Now the proof is complete.
\end{pf}

In general, the sets $U_f$ and $LU_f$ are different and the converse
of the second half of Theorem \ref{thm:char} does not hold
as the following example exhibits.

\begin{example}\label{ex:1}
Let $f(z)=(e^{\pi z}-1)/\pi.$
Then it is easy to check that $f\in\es.$
A simple computation gives us that $f(z)/zf'(z)=(e^{-\pi z}-1)/(-\pi z),$
which is obviously bounded on $\D.$
However, $f_c(z)=K_c[f](z)=z^{1-c}(e^{\pi z}-1)^c$ 
is not univalent for each $c>1.$
Indeed, since $\arg f_c(z)=(1-c)\arg z+c\arg(e^{\pi z}-1)$ for $c>0,$
one has
$$
\arg f_c(i)=(1-c)\frac\pi2+c\pi=\frac\pi2(1+c),
$$
where $i=\sqrt{-1}.$
In particular, $\arg f_c(i)>\pi/2$ for $c>1.$
This implies that $f_c(\D^+)$ intersects the negative real axis $(-\infty,0)$
for $c>1.$ Here, $\D^+=\{z\in\D: \Im z>0\}.$
Take a point $z_c\in\D^+$ so that $f_c(z_c)\in(-\infty,0)$ for $c>1.$
In view of the symmetric property that $\overline{f_c(\bar z)}=f_c(z),$
one has $f_c(z_c)=f_c(\overline{z_c}),$ which implies
$c\notin U_f$ for $c>1.$
Therefore, we conclude that $1\in U_f\setminus\Int U_f.$
\end{example}

\section{Proof of Theorem \ref{thm:qc}}\label{sec:qc}

In recent years, holomorphic motions are intensively studied in
various contexts and found a number of deep results and interesting
applications.
Among them, we need the following assertions,
which can be found, for instance, in a paper \cite{AM03} by Astala and Martin.

\begin{lem}
Let $F:D\times E\to \sphere$ be a holomorphic motion of $E\subset\sphere$
over a simply connected domain $D\subset\C$ with basepoint $c_0.$
Then $F$ extends to a holomorphic motion $\hat F$ of $\sphere$ over $D.$
Moreover, $\hat F$ is jointly continuous on $D\times\sphere$ and
$\hat F(c,\cdot)$ is $\tanh d_D(c_0,c)$-quasiconformal automorphism of $\sphere$
for each $c\in D.$
\end{lem}

We are in a position to prove Theorem \ref{thm:qc}.

\begin{pf}[Proof of Theorem \ref{thm:qc}]
Recall that $D$ is a connected component of $\Int U_f$ for a non-identity
function $f\in\ZF.$
Let $f_c=K_c[f]$ for $c\in D.$
By Theorem \ref{thm:Uf}, $D$ is simply connected.
Fix $c_0\in D$ and consider the function
$F(c,w)=(f_c\circ f_{c_0}\inv)(w)$ for $c\in D$ and $w\in f_{c_0}(\D)$
and $F(c,\infty)=\infty$ for $c\in D.$
Then this gives a holomorphic motion of $E=f_{c_0}(\D)\cup\{\infty\}$
over $D$ with basepoint $c_0.$
We now apply the above lemma to obtain an extension $\hat F$ of $F$
to $D\times\sphere.$
Then $g=\hat F(c,\cdot)$ gives a $\tanh d_D(c_0,c)$-quasiconformal automorphism
of $\C$ such that $f_c=g\circ f_{c_0}$ on $\D.$
Thus the proof is complete.
\end{pf}

\section{Proof of Theorem \ref{thm:two} and concluding remarks}\label{sec:two}

In order to construct such an example as in Theorem \ref{thm:two},
we will make use of the following
univalence criterion, which may be of independent interest.

\begin{lem}\label{lem:univ}
There exists a positive number $m$
such that the condition $e^{-m}<|zf'(z)/f(z)|<e^{m}$ on $|z|<1$
implies univalence of $f$ on $\D$ for $f\in\A.$
\end{lem}

\begin{pf}
Fix an arbitrary positive number $m$ and let $\alpha=2im/\pi.$
Then the function
$$
q(z)=\left(\frac{1+z}{1-z}\right)^\alpha
$$
is a universal covering projection of $\D$ onto the annulus
$e^{-m}<|w|<e^m$ with $q(0)=1.$
Therefore, the assumption means that the function $p(z)=zf'(z)/f(z)$ 
is subordinate to $q(z);$ in other words, 
$p=q\circ\omega$ for a function $\omega\in\A$ with $|\omega(z)|\le |z|.$
The situation is same as in the proof of Theorem 1.1 in \cite{SugawaST}
except for the exponent $\alpha,$ which is assumed to be a positive number there.
Thus we follow the argument in \cite{SugawaST}.

We first observe that the inequality
$$
\left|\log\frac{1+z}{1-z}\right|
=\left|\sum_{n=1}^\infty\frac{2z^{2n-1}}{2n-1}\right|
\le\sum_{n=1}^\infty\frac{2|z|^{2n-1}}{2n-1}
=\log\frac{1+|z|}{1-|z|}
$$
holds.
Hence, letting $w=\log((1+z)/(1-z)), W=\log((1+|z|)/(1-|z|)),$
we have
$$
|q(z)-1|=|e^{\alpha w}-1|
=\left|\sum_{n=1}^\infty\frac{(\alpha w)^n}{n!}\right|
\le\sum_{n=1}^\infty\frac{(|\alpha| W)^n}{n!}
=e^{|\alpha|W}-1=Q(|z|)-1,
$$
where $Q(z)=((1+z)/(1-z))^{|\alpha|}.$
By using this, we have the inequality
$|f''(z)/f'(z)|\le F''(|z|)/F'(|z|)$ for $|z|<1$
in the same way as in \cite{SugawaST},
where $F\in\A_1$ is defined by the relation $zF'(z)/F(z)=Q(z).$
In particular, we have
$$
\sup_{z\in\D}(1-|z|^2)\left|\frac{f''(z)}{f'(z)}\right|
\le
\sup_{z\in\D}(1-|z|^2)\left|\frac{F''(z)}{F'(z)}\right|.
$$
The right-hand term is estimated by $6|\alpha|$ from above
(cf.~\cite{SugawaST}).
Therefore, when $m\le\pi/12,$ we have
$$
\sup_{z\in\D}(1-|z|^2)\left|\frac{f''(z)}{f'(z)}\right|
\le 6|\alpha|=\frac{12}\pi m\le 1.
$$
Becker's theorem \cite{Becker72} now implies univalence
of $f.$
Thus the lemma has been proved with the choice $m=\pi/12.$
\end{pf}

We made a crude estimate above.
Therefore, $\pi/12$ is not the sharp constant.
As we will see below, $m$ cannot be taken so that $m>\pi/2.$
It may be an interesting problem to find (or to estimate)
the best possible value of $m$ in the lemma.
Since the problem is out of our scope in this note,
we will treat this problem in a separate paper.

We now prove Theorem \ref{thm:two}.

\begin{pf}[Proof of Theorem \ref{thm:two}]
Let $m$ be the number which appears in Lemma \ref{lem:univ}.
Let $f\in\A_1$ be the function determined by the relation
$zf'(z)/f(z)=((1+z)/(1-z))^{im/\pi}.$
Then
$$
V(f)=\{w\in\C: e^{-m/2}<|w|<e^{m/2}\}.
$$
Since $V(f)$ separates $0$ from $\infty,$
the image $T(V(f))$ under the M\"obius transformation $T$
given in \eqref{eq:T} separates $0$ and $1.$
Since $T(V(f))=\sphere\setminus LU_f\subset\sphere\setminus U_f$
by \eqref{eq:fund}, it is enough to see that $\{0,1\}\subset\Int U_f$
to obtain a desired example.
The first assertion of Theorem \ref{thm:char} implies
$0\in\Int U_f$ because $zf'(z)/f(z)$ is bounded.
On the other hand, in view of \eqref{eq:fc}, the range of
$zf_c'(z)/f_c(z)$ stays in the annulus $e^{-m}<|w|<e^m$ for
$c$ close enough to $1,$ where $f_c=K_c[f]$ for $c\in\C.$
Therefore, Lemma \ref{lem:univ} implies that $f_c$ is univalent on $\D$
when $|c-1|$ is small enough.
This means that $1\in\Int U_f.$
Now the program of the proof is complete.
\end{pf}

The same technique used in the above proof yields the following.

\begin{prop}
Let $A$ be a compact, polynomially convex subset of $\C$ with $0\in A.$
Then there exists a function $f$ in $\ZF$ such that
$LU_f=A.$
\end{prop}

\begin{pf}
Let $\Omega=T\inv(\sphere\setminus A).$
Then the polynomial convexity of $A$ implies that $\Omega$ is a
domain (a connected non-empty open set) in $\sphere.$
Note also that $1\in\Omega\subset\C\setminus\{0\}.$
When $A$ consists of $0$ only; namely, $\Omega=\C\setminus\{0\},$
we choose $f\in\ZF$ so that
$zf'(z)/f(z)=((1+z)/(1-z))^3$ for example.
Then $V(f)=\C\setminus\{0\}$ and thus $LU_f=A=\{0\}$ by \eqref{eq:fund}.
We now assume that $A$ contains at least two points.
Then, thanks to the uniformization theorem,
we can take a holomorphic universal covering projection 
$p$ of $\D$ onto $\Omega$ with $p(0)=1.$
If we take $f\in\ZF$ so that
$zf'(z)/f(z)=p(z),$ we have $V(f)=\Omega.$
Therefore, $LU_f=\C\setminus T(V(f))=A$ by \eqref{eq:fund}.
\end{pf}


We mention necessary conditions for univalence of $f$ in terms of
its power deformations.
Prawitz \cite{Pra27} extended Gronwall's area theorem in the following way
(see also \cite{Milin:univ}).
Let $F(\zeta)=\zeta+b_0+b_1/\zeta+b_2/\zeta^2+\dots$ be a non-vanishing
univalent meromorphic function in $|\zeta|>1.$
When
$$
\left(\frac{F(\zeta)}{\zeta}\right)^\lambda
=\sum_{n=0}^\infty D_n(\lambda)\zeta^{-n},\quad |\zeta|>1,
$$
the inequality
$$
\sum_{n=0}^\infty (\lambda-n)|D_n(\lambda)|^2
=\lambda+\sum_{n=1}^\infty (\lambda-n)|D_n(\lambda)|^2\ge0
$$
holds for each $\lambda>0.$
This result can be translated into our setting in a simple way.

\begin{thm}[A variant of Prawitz's area theorem]
Let $f\in\es$
and $K_{c}[f](z)=z+\sum_{n=2}^\infty a_n(c)z^n$ for $c\in\C.$
Then the inequality
$$
\sum_{n=2}^\infty(n-1-\lambda)|a_n(-\lambda)|^2\le \lambda
$$
holds for each $\lambda>0.$
\end{thm}

\begin{pf}
Let $F(\zeta)=1/f(1/\zeta).$
We note that 
$$
\left(\frac{F(\zeta)}{\zeta}\right)^\lambda
=\left(\frac {f(1/\zeta)}{1/\zeta}\right)^{-\lambda}
=1+\sum_{n=1}^\infty a_{n+1}(-\lambda)\zeta^{-n}.
$$
Prawitz's area theorem now yields the required inequality.
\end{pf}

We note that the coefficient $a_n(c)$ is a polynomial in $c$
for each $n.$
(This is true for a general $f\in\ZF.$)
For instance, when $zf'(z)/f(z)=((1+z)/(1-z))^\alpha,~\alpha=2im/\pi,$
we have $f_c(z)=z+2c\alpha z^2+c(1+2c)\alpha^2z^3+\dots.$
We now suppose that $f$ is univalent in $\D.$
Taking only $a_2(-\lambda)$-term in the last theorem, 
we obtain the inequality
$$
4\lambda^2(1-\lambda)|\alpha|^2\le \lambda,
$$
which is equivalent to $4\lambda(1-\lambda)|\alpha|^2\le 1.$
Letting $\lambda=1/2,$ we have $|\alpha|\le1.$
In this way, we showed that $m\le \pi/2$ is necessary for univalence
of $f.$
We could improve this upper bound if we increase the number of
terms in the summation with the cost of calculation amount.

\def\cprime{$'$} \def\cprime{$'$} \def\cprime{$'$}
\providecommand{\bysame}{\leavevmode\hbox to3em{\hrulefill}\thinspace}
\providecommand{\MR}{\relax\ifhmode\unskip\space\fi MR }
\providecommand{\MRhref}[2]{%
  \href{http://www.ams.org/mathscinet-getitem?mr=#1}{#2}
}
\providecommand{\href}[2]{#2}

\end{document}